\documentclass[a4paper,10pt]{article}
\usepackage{stmaryrd}
\usepackage{amsfonts}
\usepackage{bbm}
\usepackage{amscd}
\usepackage{mathrsfs}
\usepackage{latexsym,amssymb,amsmath,amscd,amscd,amsthm,amsxtra,xypic}
\usepackage[dvips]{graphicx}
\usepackage[utf8]{inputenc}
\usepackage[T1]{fontenc}
\usepackage{lmodern}
\usepackage{amssymb}
\usepackage[all]{xy}
\usepackage{nicefrac,mathtools,enumitem}
\usepackage{microtype}

\newtheorem{thm}{Theorem}[section]
\newtheorem{lem}[thm]{Lemma}
\newtheorem{cor}[thm]{Corollary}
\newtheorem{pro}[thm]{Proposition}
\newtheorem{ex}[thm]{Example}
\newtheorem{ep}[thm]{Example}
\newtheorem{rmk}[thm]{Remark}
\newtheorem{defi}[thm]{Definition}

\newcommand{\be }{\begin{equation}}
\newcommand{\ee }{\end{equation}}

\newcommand{\pf}{\noindent{\bf Proof.}\ }

\newcommand{\g}{\mathbbm g}
\newcommand{\h}{\mathbbm h}
\newcommand{\hh}{\mathrm h}
\newcommand{\G}{\mathbb G}
\newcommand{\HH}{\mathbb H}

\newcommand{\R}{\mathbb R}

\newcommand{\huaG}{\mathcal{G}}
\renewcommand{\k}{\mathbbm k}

\newcommand{\frke}{\mathfrak e}

\newcommand{\frkg}{\mathfrak g}
\newcommand{\frkh}{\mathfrak h}

\newcommand{\frkk}{\mathfrak k}

\newcommand{\frkG}{\mathfrak G}

\def\qed{\hfill ~\vrule height6pt width6pt depth0pt}

\newcommand{\half}{\frac{1}{2}}

\newcommand{\pair}[1]{\left\langle #1\right\rangle}

\newcommand{\Courant}[1]{\left\llbracket  #1\right\rrbracket }

\newcommand{\Id}{{\rm{Id}}}

\newcommand{\id}{\mathbbm{i}}

\newcommand{\dM}{\mathrm{d}}

\newcommand{\Der}{\mathrm{Der}}

\newcommand{\Aut}{\mathrm{Aut}}

\newcommand{\gl}{\mathfrak {gl}}

\newcommand{\Ker}{\mathrm{Ker}}

\newcommand{\End}{\mathrm{End}}
\newcommand{\ad}{\mathrm{ad}}
\newcommand{\pr}{\mathrm{pr}}

\newcommand{\ve}{\mathrm{v}}

\newcommand{\sgn}{\mathrm{sgn}}
\newcommand{\Ksgn}{\mathrm{Ksgn}}

\newcommand{\V}{\mathbb{V}}

\begin{document}
\title{
{Integration of semidirect product Lie 2-algebras
\thanks
 {The first author is supported
by NSFC  (11026046,11101179), SRFDP (20100061120096)
 and ``the Fundamental
Research Funds for the Central Universitie'' (200903294).  The
second author is supported by the German Research Foundation
(Deutsche Forschungsgemeinschaft (DFG)) through the Institutional
Strategy of the University of G\"ottingen.
 }
} }
\author{Yunhe Sheng  \\
Department of Mathematics, Jilin University,
 Changchun 130012, Jilin, China
\\\vspace{3mm}
email: shengyh@jlu.edu.cn\\
Chenchang Zhu\\
Courant Research Center ``Higher Order Structures'', University of
G$\ddot{\rm{o}}$ttingen\\
email:zhu@uni-math.gwdg.de}
\date{}
\footnotetext{{\it{Keyword}:  representations up to homotopy,
$L_\infty$-algebras, Lie 2-algebras, omni-Lie algebras, crossed
modules, integration, butterfly}}

\footnotetext{{\it{MSC}}: Primary 17B55. Secondary 18B40, 18D10.}

\maketitle
\begin{abstract}
The semidirect product of a Lie algebra and a 2-term representation up to
homotopy is a Lie 2-algebra. Such Lie 2-algebras include many examples arising from the Courant algebroid
appearing in generalized complex geometry. In this paper, we integrate  such a Lie 2-algebra to a strict Lie 2-group in the finite dimensional case.
\end{abstract}
\tableofcontents
\section{Introduction}

In recent years people have paid much attention to the integration
of  Lie-algebra-like structures, such as that of Lie algebroids
\cite{cf, tz2}, of $L_\infty$-algebras \cite{getzler, henriques} and
of Courant algebroids \cite{LS,MT,shengzhu3}. Here ``integration''
is meant in the same sense in which a Lie algebra is integrated to a
corresponding Lie group. However, directly applying the general
method proposed in \cite{getzler, henriques} to certain examples,
such as Courant algebroids and its linearization---omni-Lie
algebras, gives an infinite dimensional space which is too abstract
for our application. In this paper, we give an explicit finite
dimensional Lie 2-group integrating   a certain type of Lie
2-algebras in the form of semidirect products. We are motivated to
integrate such a type of Lie 2-algebras because it contains many
meaningful examples, such as Courant algebroids, omni-Lie
algebroids/omni-Lie algebras, string Lie 2-algebras, as shown in
\cite{shengzhu1}. We give an integration of omni-Lie algebras at the
end of the paper. The integration of Courant algebroids and string
Lie 2-algebras is given in a later paper \cite{shengzhu3}.

Just as we can form a semidirect product Lie algebra $\frkg \ltimes
V$ out of  a representation of a Lie algebra, we can form a new
$L_\infty$-algebra via the semidirect product construction. For
this, we need to replace a usual representation by an
$L_\infty$-module in the sense of \cite{lada-markl}. Notice that
even for a Lie algebra $\frkg$, its $L_\infty$-module is a more
general concept than its usual representation. These
$L_\infty$-modules are also studied by
\cite{abad-crainic:rep-homotopy} under the name ``representation up
to homotopy'' in the context of Lie algebroids and Lie groupoids.
The semidirect product of a Lie algebra $\frkg$ with
its 2-term $L_\infty$-module gives us a Lie 2-algebra. It is this
sort of Lie 2-algebras we integrate in this paper.

 It is rather easy to integrate a strict morphism between
strict Lie 2-algebras. But to integrate an $L_{\infty}$-morphism
between strict Lie 2-algebras is not that straightforward. We have
to solve a set of PDE's, which is a modification of the one that Lie
solved to integrate Lie algebra morphisms. We have succeeded in
solving them \cite{shengzhu4}, but the model we obtained is infinite
dimensional. However, thanks to Noohi's timely paper on butterfly
method of $L_{\infty}$-morphisms \cite{noohi:morphism}, where he
realizes an $L_\infty$-morphism of strict Lie 2-algebras as a
zig-zag of strict morphisms, we are able to avoid the
above-mentioned technical difficulties.

Here is our integration procedure: first we form the semidirect product in the strict case both at the
Lie group level and at the Lie algebra level (Corollary \ref{cor:cm
group} and Corollary \ref{thm:semidirect product of strict 1}).
Using the butterfly, we strictify an $L_\infty$-morphism between
strict Lie 2-algebras. Then by using the strict morphism, we form a
semidirect product Lie 2-algebra, which is equivalent to the
original one (Theorem \ref{thm:equivalent}). Thus we obtain a strict
Lie 2-group integrating the semidirect product at the Lie algebra
level (Corollary \ref{thm:integration}). Using the zig-zag provided
by butterfly methods, we eventually solve our integration problem
within the finite dimensional world.

As an example, we apply our integration to omni-Lie algebras.

Finally, since we realize the integration of a semidirect product Lie
2-algebra  as a strict Lie 2-group, we also provide a
strictification of the Lie 2-group we constructed in
\cite{shengzhu1} integrating the string type Lie 2-algebra $\R \to
\frkg\oplus \frkg^*$. Note that there is not yet a canonical way
to  strictify a Lie 2-group, so we are quite lucky to achieve this.

The paper is organized as follows. In Section 2 we briefly recall
some notions and basic facts about $L_\infty$-algebras,
representation up to homotopy of Lie algebras, crossed modules of
Lie algebras (resp. Lie groups), Lie 2-groups, and butterflies. In
Section 3 first we give the integration of the crossed module of Lie
algebras $\End(\V)$, which turns out to be a crossed module of Lie
groups $\Aut(\V)$. Then we construct a strict Lie 2-group associated
to any strict morphism  of crossed modules of Lie groups from
$(H_1,H_0,t,\Phi)$ to $\Aut(\V)$ (Theorem \ref{thm:semidirect
product of strict}). Then with the help of the butterfly, we give
the integration of semidirect product Lie 2-algebras. In Section 4
we give the integration of omni-Lie algebras as an application.

  {\bf
Notations:} $\dM$ is the differential in a complex of vector spaces,
$\delta$ is the differential in the DGLA associated with a complex
of vector spaces. $i$ is the inclusion map. $i_1$  and $i_2$ are the
inclusion to the first factor and the second factor respectively.
$\Id$ is the identity map. For a graded vector space $V_\bullet=
\sum_n V_n $, $V[l]$ denote the $l$-shifted graded vector space,
namely $V[l]_n=V_{l+n}$; $Sym(V)$ is the symmetric power of $V$.
$\V$ is a 2-term complex $V_1\stackrel{\dM}{\longrightarrow}V_0.$

{\bf Acknowledgement:} We give warmest thanks to John Baez, Joel
Kamnitzer, Behrang Noohi,  Weiwei Pan, Urs Schreiber, Giorgio
Trentinaglia, and Alan Weinstein for very useful comments. Y. Sheng
gratefully acknowledges the support of Courant Research Center
``Higher Order Structures'', G\"{o}ttingen University, where this
work was done during his visit.

\section{Preliminaries}
\subsection{Lie 2-algebras}
\begin{defi}\label{def:Linfinity}
An $L_\infty$-algebra is a graded  vector space $L=L_0\oplus
L_1\oplus\cdots$ equipped with a system $\{l_k|~1\leq k<\infty\}$ of
linear maps $l_k:\wedge^kL\longrightarrow L$ of degree
$\deg(l_k)=k-2$, where the exterior powers are interpreted in the
graded sense and the following relation with Koszul sign ``Ksgn'' is
satisfied for all $n\geq0$:
\begin{equation} \label{eq:jacobi}
\sum_{i+j=n+1}(-1)^{i(j-1)}\sum_{\sigma}\sgn(\sigma)\Ksgn(\sigma)l_j(l_i(x_{\sigma(1)},\cdots,x_{\sigma(i)}),x_{\sigma(i+1)},\cdots,x_{\sigma(n)})=0,
\end{equation}
where the summation is taken over all $(i,n-i)$-unshuffles with
$i\geq1$.
\end{defi}

For $n=1$, we have
$$
l_1^2=0,\quad l_1:L_{i+1}\longrightarrow L_i,
$$
which means that $L$ is a complex, so  we  write $\dM=l_1$ as usual.
For $n=2$, we have
$$
\dM l_2(X,Y)=l_2(\dM X,Y)+(-1)^pl_2(X,\dM Y),\quad \forall~X\in L_p,
Y\in L_q,
$$
which means that $\dM $ is a derivation with respect to $l_2$. We
view $l_2$ as a bracket $[\cdot,\cdot]$. However, it is not
a Lie bracket, the obstruction of Jacobi identity is controlled by
$l_3$ and $l_3$ also satisfies higher coherence laws.

In particular, if the $k$-ary brackets are zero for all $k>2$, we
recover the usual notion of {\bf differential graded Lie algebras}
(DGLA). If $L$ is concentrated in degrees $<n$, we get the notion of
{\bf $n$-term $L_\infty$-algebras}. A $2$-term $L_\infty$-algebra is
also called a {\bf Lie 2-algebra} in this paper. A 2-term DGLA is a special Lie 2-algebra with $l_3=0$, thus it is also
called a {\bf strict Lie 2-algebra}.  For details about
Lie 2-algebras, see
\cite{baez:2algebras,baez:classicalstring,roytenbergweakLie2}.

\begin{defi}\label{def:m of 2-alg}
Let $L$ and $L'$ be Lie 2-algebras. A Lie 2-algebra
morphism\footnote{Since we view 2-term $L_\infty$-algebras as Lie
2-algebras, our Lie 2-algebra morphisms are exactly
$L_\infty$-algebra morphisms.} $f:L\longrightarrow L'$ consists of:
\begin{itemize}
\item[$\bullet$]two linear maps $f_0:L_0\longrightarrow L'_0$ and
$f_1:L_1\longrightarrow L'_1$ preserving the differential;
\item[$\bullet$]a skew-symmetric bilinear map
$f_2:\wedge^2L_0\longrightarrow L'_1$,
\end{itemize} such that the following equalities hold for all $X,Y,Z\in L_0,~A\in
L_1$,
\begin{eqnarray}\label{eqn:DGLA morphism c 1}
f_0(l_2(X,Y))-l_2'(f_0(X),f_0(Y))&=&\dM' f_2(X,Y),\\\nonumber
f_1(l_2(X,A))-l'_2(f_0(X),f_1(A))&=&f_2(X,\dM A),\\\nonumber
l'_2(f_0(X),f_2(Y,Z))+c.p.+f_1(l_3(X,Y,Z))&=&f_2(l_2(X,Y),Z)+c.p.+l_3'(f_0(X),f_0(Y),f_0(Z)),
\end{eqnarray}
where $c.p.$ means cyclic permutations. A morphism $f$ is called
{\bf strict} if $f_2=0$.
\end{defi}

In particular, if $L$ and $W$ are strict Lie 2-algebras, we can
still have non-strict morphisms between them. Only the last equality
in the above definition simplifies to
\begin{equation}\label{eqn:DGLA morphism c 2}
l_2'(f_0(X),f_2(Y,Z))+c.p.=f_2(l_2(X,Y),Z)+c.p..
\end{equation}

By a complex of vector spaces, we mean a graded vector space
$V_\bullet$ endowed with a degree minus one endomorphism $\dM$
satisfying $\dM^2=0$:
\begin{equation}
(V_\bullet,\dM):\cdots\stackrel{\dM}{\longrightarrow}V_k\stackrel{\dM}{\longrightarrow}V_{k-1}
\stackrel{\dM}{\longrightarrow}\cdots\cdots\stackrel{\dM}{\longrightarrow}V_0.
\end{equation}
An element $u\in V_k$ is called a homogeneous element of degree $k$.
From this graded vector space, we can form a new graded vector space
$\End(V_\bullet)$, the degree $k$ part $\End^k(V_\bullet)$ of which
consists of linear maps $T:V_\bullet\longrightarrow V_\bullet$ which
increase the degree by $k$. We denote $|T|$ the degree of $T$. We
introduce an operator of  degree minus one $\delta$ on
$\End(V_\bullet)$ by setting:
\begin{equation}\label{eqn:differential}
\delta(T)=\dM\circ T-(-1)^{k}T\circ\dM,\quad\forall
~T\in\End^k(V_\bullet).
\end{equation}
We have also the super-commutator bracket $[\cdot,\cdot]$ on
$\End(V_\bullet)$ given by the linear expansion of the formula for
homogeneous elements:
\begin{equation}\label{eqn:super bracket}
[T,S]=T\circ S-(-1)^{|T||S|}S\circ T.
\end{equation}

\begin{pro}\label{thm:dgla}
With the above notations,  $(\End(V_\bullet),[\cdot,\cdot],\delta)$
is a DGLA, where $\delta$ and $[\cdot,\cdot]$ are given by
$(\ref{eqn:differential})$ and $(\ref{eqn:super bracket})$
respectively.
\end{pro}

From now on, we focus on a complex of vector spaces of length 2,
$\V: V_1\stackrel{\dM}{\longrightarrow}V_0$. By
 Proposition \ref{thm:dgla}, we obtain a DGLA of length 3:
\begin{equation}\label{eqn:dgla of 3}
\End(V_0,V_1)\stackrel{\delta}{\longrightarrow}\End(V_0,V_0)\oplus\End(V_1,V_1)\stackrel{\delta}{\longrightarrow}\End(V_1,V_0),
\end{equation}
where the degree-1 part $\End^1(\V)$ is $\End(V_0,V_1)$, the
degree-$0$ part $\End^0(\V)$ is $\End(V_0,V_0)\oplus\End(V_1,V_1)$
and the degree-$(-1)$ part $\End^{-1}(\V)$ is $\End(V_1,V_0)$. Since
$\delta^2=0$, we have
$$\delta(\End^1(\V))\subset\Ker(\delta|_{\End^0(\V)}),$$
Furthermore, since $\delta$ is a derivation,
$\Ker(\delta|_{\End^0(\V)})$ is a Lie subalgebra of $\End^0(\V)$. In
fact, by (\ref{eqn:differential}), the definition of $\delta$, we
have
$$
\Ker(\delta|_{\End^0(\V)})=\{(A,B)\in\End(V_0,V_0)\oplus\End(V_1,V_1)|A\circ\dM=\dM\circ
B\}.
$$
We denote $\Ker(\delta|_{\End^0(\V)})$ by $\End^0_\dM(\V)$. Then, by
truncation, we obtain a new DGLA of length 2, which we denote by
$\End(\V)$,
\begin{equation}\label{eqn:dgla of 2}
\End(\V):\End^1(\V)\stackrel{\delta}{\longrightarrow}\End^0_\dM(\V).
\end{equation}

It is well known that an $L_\infty$-algebra structure\footnote{This works for finite dimensional $L_i$'s which is
our setting in this paper.} on $L$ is equivalent
to a differential graded commutative algebra (d.g.c.a) structure on
$$ Sym (L^*[-1])= \underbrace{\k}_{{\text{degree}\; 0}} \oplus
\underbrace{L^*_0}_{{\text{degree}\; 1}} \oplus \underbrace{\big[
\wedge^2 L^*_0 \oplus L^*_{-1}\big]}_{\text{degree} \;2} \oplus
\underbrace{ \big[ \wedge^3 L^*_0 \oplus L^*_0 \otimes L^*_{-1}
\oplus L^*_{-2}\big]}_{\text{degree} \;3} \oplus \dots.
$$
The degree $1$ differential $\delta$ is given by dualizing
$l_i's$. An {\bf $L_\infty$-module}  (or a {\bf representation up to
  homotopy} \cite{abad-crainic:rep-homotopy}) of an $L_\infty$-algebra $L$ is a graded vector space
$V_\bullet$ with a degree one differential $D$ on the graded
vector space
$$ (Sym (L^*[-1]) \otimes V_\bullet)_n = \oplus_k Sym(L^*[-1])_k \otimes
V_{n-k}. $$ We denote the space of $L_\infty$-modules of $L$ by
$Rep^\infty(L)$. If $V_\bullet \in Rep^\infty(L)$, then the  dual
$V_\bullet^* $ and symmetric power $Sym (V_\bullet^*)$  are in
$Rep^\infty(L)$ too. By the fact that
$$
Sym(L^*[-1])\otimes Sym(V_\bullet^*[-1])\cong Sym\big((L\oplus
V_\bullet)^*[-1]\big),
$$
we obtain an $L_\infty$-algebra structure on $L\oplus V_\bullet$
from the $L_\infty$-module structure on $V_\bullet$. This natural
$L_\infty$-algebra structure on $L\oplus V_\bullet$  is called
the {\bf semidirect product} of $L$ and
$V_\bullet$.

An {\bf equivalent definition of $L_\infty$-modules} is given
in \cite{lada-markl} using $L_\infty$-morphisms.  In the following,
we focus on the case that $L$ is simply a Lie algebra or a  strict Lie
2-algebra and
$V_\bullet$ is a 2-term complex of vector spaces. In this case, we
use the terminology ``representation up to homotopy'' because this
terminology can also apply to the situation of groups or groupoids,
which is the case we will use later.

Thus explicitly in our situation, a 2-term representation up to
homotopy of a Lie algebra $\frkg$ on a 2-term complex of vector
spaces $\V$ is a Lie 2-algebra morphism $(\mu,\nu)$ from the Lie
algebra $\frkg$ to the $2$-term DGLA $\End(\V)$,
$$
 \xymatrix{
\nu:\wedge^2\frkg\ar[r]&\End^1(\V)\ar[d]_{\delta}\\
 \mu:\frkg\ar[r]& \End^0_\dM(\V).}
$$
 Then the corresponding
semidirect product Lie 2-algebras are given explicitly as follows.
The 2-term complex of vector spaces is given by:
$$
\frkg\ltimes \V:V_1\stackrel{0+\dM}{\longrightarrow}\frkg\oplus V_0;
$$
$l_2:\wedge^2(\frkg\ltimes \V)\longrightarrow\frkg\ltimes \V $ is
given by:
\begin{equation}\label{2bracket}
\left\{\begin{array}{rll}l_2((X,\xi),(Y,\eta))&=&\big([X,Y],\mu_0(X)(\eta)-\mu_0(Y)(\xi)\big),\\
l_2((X,\xi),m)&=&\mu_1(X)(m),\\
l_2(m,n)&=&0,\end{array}\right.
\end{equation}
for any $(X,\xi),(Y,\eta)\in\frkg\oplus V_0,~m,~n\in V_1$, and
$l_3(\cdot,\cdot,\cdot):\wedge^3(\frkg\ltimes
\V)\longrightarrow\frkg\ltimes \V $ is given by:
\begin{equation}\label{3bracket}
l_3((X,\xi),(Y,\eta),(Z,\gamma))=-\nu(X,Y)(\gamma)+c.p..
\end{equation}
One should note that $l_2(\cdot,\cdot)$ is not a Lie bracket.
Instead, we have
$$
l_2(l_2((X,\xi),(Y,\eta)),(Z,\gamma))+c.p.=\dM(\nu(X,Y)(\gamma))+c.p..
$$
\subsection{2-term DGLAs and crossed modules of Lie algebras}

\begin{defi}
A crossed module of Lie algebras is a quadruple
$(\frkh_1,\frkh_0,dt,\phi)$, which we denote by $\h$, where
$\frkh_1$ and $\frkh_0$ are Lie algebras,
$dt:\frkh_1\longrightarrow\frkh_0$ is a Lie algebra morphism and
$\phi:\frkh_0\longrightarrow\Der(\frkh_1)$ is an action of Lie
algebra $\frkh_0$ on Lie algebra $\frkh_1$ as a derivation, such
that
$$
dt(\phi_X(A))=[X,dt(A)],\quad \phi_{dt(A)}(B)=[A,B].
$$
\end{defi}

\begin{ex}
For any Lie algebra $\frkk$, the adjoint action $\ad$ is a Lie
algebra morphism from $\frkk$ to $\Der(\frkk)$. Then
$(\frkk,\Der(\frkk),\ad,\Id)$ is a crossed module of Lie algebras.
\end{ex}

\begin{thm}\label{thm:dgla and cm}
There is a one-to-one correspondence between 2-term DGLAs and
crossed modules of Lie algebras.
\end{thm}

In short, the formula for the correspondence can be given as follows: A 2-term DGLA $L_1\stackrel{\dM}{\longrightarrow} L_0$
gives rise to   a Lie algebra crossed module with $\frkh_1=L_1$ and
$\frkh_0=L_0$, where the Lie brackets are given by:
\begin{eqnarray*}
~[A,B]_{\frkh_1}&=&l_2(\dM(A),B),\quad\forall~A,B\in
L_1,\\
~[X,X^\prime]_{\frkh_0}&=&l_2(X,X^\prime),\quad\forall~X,X^\prime\in
L_0,
\end{eqnarray*}
and $dt=\dM$, $\phi:\frkh_0\longrightarrow\Der(\frkh_1)$ is given by
$ \phi_X(A)=[X,A]. $ The DGLA structure gives the Jacobi identity for
$[\cdot,\cdot]_{\frkh_1}$ and $[\cdot,\cdot]_{\frkh_0}$, and various other
conditions for crossed modules.

Conversely, a crossed module $(\frkh_1,\frkh_0,dt,\phi)$ gives rise
to a 2-term DGLA with $\dM=dt$, $L_1=\frkh_1$ and $L_0=\frkh_0$, and
$l_2's$ are given by:
\begin{eqnarray*}
~l_2(A,B)&\triangleq& 0,~\qquad\qquad \forall ~A,B\in\frkh_1,\\
~ l_2(X,X^\prime)&\triangleq&[X,X^\prime]_{\frkh_0},\quad \forall
~X,X^\prime\in\frkh_0,\\
~l_2(X,A)&\triangleq&\phi_X(A).
\end{eqnarray*}

The corresponding crossed module of the 2-term DGLA $\End(\V)$ given
by (\ref{eqn:dgla of 2}) is as follows: the Lie algebra $\frkk_1$ as
a vector space is $\End^1(\V)$. Its Lie bracket is given by
\begin{equation}\label{k1}
[A,B]_{\frkk_1}=[\delta(A),B]=A\circ\dM\circ B-B\circ\dM\circ A.
\end{equation} The Lie algebra $\frkk_0$ is
the Lie sub-algebra $\End^0_{\dM}(\V)$ of $\End^0(\V)$,
\begin{equation}\label{k0}
\frkk_0\triangleq\End^0_\dM(\V),\quad [X,X^\prime]_{\frkk_0}=X\circ
X^\prime-X^\prime\circ X.
\end{equation}
Furthermore, the Lie algebra morphism  $dt$ and the action $\phi$
are given by
\begin{equation} \label{eq:phi}
\phi_X(A)=[X,A],\quad dt=\delta.
\end{equation}
We denote this crossed module of Lie algebras also by $\End(\V)$,
i.e.
\begin{equation}\label{crossed module endv}
\End(\V)=(\frkk_1,\frkk_0,dt,\phi),
\end{equation}
where $\frkk_1$ and $\frkk_0$ are given by (\ref{k1}) and (\ref{k0})
respectively.

\begin{ep}[Semidirect products via strict morphisms] \label{ep:semidirect}
Given a strict morphism $(\psi_1,\psi_0)$ of strict Lie 2-algebras
from the strict Lie 2-algebra corresponding  to the crossed modules
of Lie algebras $\h=(\frkh_1, \frkh_0, dt, \phi)$ to $\End(\V)$, the
semidirect product $\h \ltimes \V$ is again a strict Lie 2-algebra
whose corresponding crossed module is
\begin{equation}\label{eq:smdr-prd-alg}
(\frkh_1\ltimes V_1,\frkh_0\ltimes V_0,dt\times \dM,\phi).
\end{equation}
Here the Lie algebra structure on $\frkh_1\ltimes V_1$ is given by
$$
[(A,m),(B,n)]=\big([A,B],\psi_0(dt(A))(n)-\psi_0(dt(B))(m)\big),
$$
the Lie algebra structure on $\frkh_0\ltimes V_0$ is the semidirect
product
$$
[(X,\xi),(Y,\eta)]=\big([X,Y],\psi_0(X)(\eta)-\psi_0(Y)(\xi)\big),
$$
the Lie algebra morphism $dt\times \dM$ is given by
$$
(dt\times \dM)(A,m)=\big(dt(A),\dM m\big),
$$
and the action $\phi$ is given by
$$
\phi_{(X,\xi)}(A,m)=\big(\phi_XA,\psi_0(X)(m)-\psi_1(A)(\xi)\big).
$$
\end{ep}

\subsection{Strict Lie 2-groups and crossed modules of Lie groups}
A group is a monoid where every element has an inverse. A 2-group is
a monoidal category where every object has a weak inverse and every
morphism has an inverse. Denote the category of smooth manifolds and
smooth maps by $\rm Diff$, a (semistrict) Lie 2-group is  a 2-group
in $\rm DiffCat$, where  $\rm DiffCat$ is the 2-category consisting
of categories, functors, and natural transformations in $\rm Diff$.
For more details, see \cite{baez:2gp}. Here we only give the
definition of strict Lie 2-groups.
\begin{defi}
A strict Lie 2-group is a Lie groupoid $C$
such that
\begin{itemize}
\item[\rm(1)]
The space of morphisms $C_1$ and the space of objects $C_0$ are Lie
groups.
\item[\rm(2)] The source and the target $s,t:C_1\longrightarrow C_0$, the identity assigning function $i:C_0\longrightarrow C_1$ and
 the composition $\circ:C_1\times_{C_0}C_1\longrightarrow C_1$ are all Lie group morphisms.
\end{itemize}
\end{defi}

In the following we will denote the composition $\circ$ in the Lie
groupoid structure by $\cdot_\ve:C_1\times_{C_0}C_1\longrightarrow
C_1$ and call it the vertical multiplication. Denote the Lie 2-group
multiplication by $\cdot_\hh:C\times C\longrightarrow C$ and call it
the horizontal multiplication.

It is well known that strict Lie 2-groups can be described by
crossed modules of Lie groups.

\begin{defi}
A crossed module of Lie groups is a quadruple $(H_1,H_0,t,\Phi)$,
which we denote simply by $\HH$, where $H_1$ and $H_0$ are Lie
groups, $t:H_1\longrightarrow H_0$ is a morphism, and
$\Phi:H_0\times H_1\longrightarrow H_1 $ is an action of $H_0$ on
$H_1$ preserving the Lie group structure of $H_1$ such that the Lie
group morphism $t$ is $H_0$-equivariant:
\begin{equation}\label{cm g 1}
t\Phi_g(h)=gt(h)g^{-1},\quad \forall ~g\in H_0,~h\in H_1,
\end{equation}
and $t$ satisfies the so called Peiffer identity:
\begin{equation}\label{cm g 2}
\Phi_{t(h)}(h^\prime)=hh^\prime h^{-1},\quad\forall ~h,h^\prime\in
H_1.
\end{equation}
\end{defi}

\begin{thm}\label{thm:cm and slg}
There is a one-to-one correspondence between crossed modules of Lie
groups and strict Lie 2-groups.
\end{thm}
 Roughly speaking, given a crossed
module $(H_1,H_0,t,\Phi)$ of Lie groups, there is a strict Lie
2-group for which $C_0=H_0$ and $C_1= H_0\ltimes H_1$, the
semidirect product of $H_0$ and $H_1$. In this strict Lie 2-group,
the source and target maps $s,~t:C_1\longrightarrow C_0$ are given
by
$$
s(g,h)=g,\quad t(g,h)= t(h)\cdot g,
$$
the vertical multiplication $ \cdot_\ve$ is given by:
\begin{equation}\label{m v}
(g^\prime,h^\prime)\cdot_\ve(g,h) =(g, h^\prime\cdot h),\quad
\mbox{where} \quad g^\prime= t(h)\cdot g,
\end{equation}
the horizontal multiplication $\cdot_\mathrm{h}$ is given by
\begin{equation}\label{m h}
(g,h)\cdot_\mathrm{h} (g^\prime,h^\prime)=(g\cdot
g^\prime,h\cdot\Phi_g h^\prime).
\end{equation}

\subsection{Morphisms of crossed modules of Lie algebras and butterflies}\label{sec:butterfly}

Since 2-term DGLAs are the same as crossed modules of Lie algebras,
it is straightforward to obtain the definition of morphisms of
crossed modules of Lie algebras. Let
$$~\g=(\frkg_1,\frkg_0,dt,\phi),\quad\h=(\frkh_1,\frkh_0,dt,\phi)$$be
crossed modules of Lie algebras. Here we use the same notations $dt,
\phi$. This will not lead to confusion since the correct
interpretation will always be clear from the context.

\begin{defi}
A  morphism $f:\g\longrightarrow \h$ consists of:
\begin{itemize}
\item[$\bullet$]two linear maps $f_0:\frkg_0\longrightarrow \frkh_0$ and
$f_1:\frkg_1\longrightarrow \frkh_1$ preserving the morphism $dt$,
\item[$\bullet$]a skew-symmetric bilinear map
$f_2:\wedge^2\frkg_0\longrightarrow \frkh_1$,
\end{itemize} such that the following equalities hold for all $X,Y,Z\in \frkg_0,~A\in
\frkg_1$,
\begin{eqnarray*}
f_0[X,Y]-[f_0(X),f_0(Y)]&=&\dM f_2(X,Y),\\
f_1(\phi_XA)-\phi_{f_0(X)}f_1(A)&=&f_2(X,dt( A)),\\
~[f_0(X),f_2(Y,Z)]+c.p.&=&f_2([X,Y],Z)+c.p..
\end{eqnarray*}
\end{defi}
The morphism $f$ is called a {\em strict morphism} if $f_2=0.$

\begin{defi}
Two Lie 2-algebras are said to be equivalent if there is a Lie
2-algebra morphism which induces an equivalence of the underlying
2-term complexes of vector spaces.
\end{defi}
In particular, the {\em equivalence of crossed modules of Lie algebras} is
defined to be the equivalence of the corresponding Lie 2-algebras.\vspace{2mm}

The theory of ``butterflies'', developed by Aldrovandi and Noohi in
\cite{noohi:butterfly, noohi:morphism},  is a nice way to describe
(non-strict) morphisms of crossed modules of Lie algebras and of Lie
groups.

\begin{defi}{\rm\cite{noohi:morphism}}
A butterfly from $\g$ to $\h$ is a commutative diagram
$$
 \xymatrix{
\frkg_1\ar[dd]\ar[dr]^{\kappa}&&\frkh_1\ar[dl]_{\iota}\ar[dd]\\
&\frke\ar[dl]^{\sigma}\ar[dr]_{\rho}&\\
 \frkg_0&& \frkh_0,}
$$
in which both the diagonal sequence are complexes of Lie algebras
and the NE-SW sequence is short exact, such that for every
$A\in\frkg_1,B\in\frkh_1$ and $e\in\frke$, we have
$$
[e, \kappa(A)] = \kappa(\phi_{\sigma(e)}A),\quad [e, \iota(B)] =
\iota(\phi_{\rho(e)}B).
$$
\end{defi}

\begin{defi}{\rm\cite{noohi:butterfly}}
Let $\G$ and $\HH$ be two  crossed modules of Lie groups. A
butterfly $E:\G\longrightarrow \HH$ is a commutative diagram
\begin{equation}\label{eqn:butterfly of group}
 \xymatrix{
G_1\ar[dd]\ar[dr]^{\kappa}&&H_1\ar[dl]_{\iota}\ar[dd]\\
&E\ar[dl]^{\sigma}\ar[dr]_{\rho}&\\
 G_0&& H_0 }
\end{equation}
in which both diagonal sequences are complexes of Lie groups, and
the NE-SW sequence is short exact. Furthermore, for every $x\in
E,~\alpha\in H_1,~\beta\in G_1$, the following  equalities hold:
$$
\iota(\Phi_{\rho(x)}
\alpha)=x\iota(\alpha)x^{-1},\quad\kappa(\Phi_{\sigma(x)}
\beta)=x\kappa(\beta)x^{-1}.
$$ The butterfly $E$ is an equivalence between $\G$ and $\HH$ if and only if the NW-SE sequence is also short exact.
\end{defi}

\begin{rmk} \label{rk:morita}
 For people who understand crossed modules as Lie 2-groups, the butterfly $E$ above is an H.S. morphism between the Lie groupoids
 $G_1\times G_0 \Rightarrow G_0$ and $H_1 \times H_0 \Rightarrow
 H_0$. Moreover this H.S. morphism preserves the 2-group structure
 maps. This coincides with the notion of generalized morphisms between
 Lie 2-groups in \cite{z:tgpd-2}. The notion of equivalence coincides with the notion of Morita equivalence therein. In fact  $E$ is an
 equivalence if and only if $E$ is a Morita bibundle
 of the underling groupoids $G_1\times G_0 \Rightarrow G_0$ and $H_1 \times H_0 \Rightarrow H_0$.
\end{rmk}

It is easy to see that
\begin{cor}\label{cor:equivalence-gp}
A crossed module of Lie groups $\G=(G_1, G_0, t, \Phi)$ is
equivalent to the Lie group $G_0/G_1$ viewed as a trivial crossed
module $(1, G_0/G_1, 1, 1)$ if and only if the Lie group morphism
$t$ is injective.
\end{cor}

In \cite{noohi:butterfly}, Noohi has proved that any morphism
between crossed modules of Lie algebras can be integrated to a butterfly
of crossed modules of Lie groups. More precisely, for any morphism
$f:\g\longrightarrow\h$, where $f=(f_0,f_1,f_2),$ he defines a
bracket $[\cdot,\cdot]$ on $\frkg_0\oplus \frkh_1$,
\begin{equation}\label{bracket of B}
[(X,A),(Y,B)]=\big([X,Y],[A,B]+\phi_{f_0(X)}B-\phi_{f_0(Y)}A-f_2(X,Y)\big),
\end{equation} and  maps $\kappa,i_2,\sigma,\rho,$
\begin{eqnarray*}
\kappa:\frkg_1\longrightarrow\frkg_0\oplus
\frkh_1,\quad\kappa(A)&=&(-f_1(A),dt(A)),\\
i_2:\frkh_1\longrightarrow\frkg_0\oplus
\frkh_1,\quad i_2(B)&=&(0,B),\\
\sigma:\frkg_0\oplus
\frkh_1\longrightarrow\frkg_0,\quad\sigma(X,B)&=&X,\\
\rho:\frkg_0\oplus
\frkh_1\longrightarrow\frkh_0,\quad\rho(X,B)&=&f_0(X)+dt(B).
\end{eqnarray*}
Then he  obtains a butterfly
\begin{equation}\label{eqn:butterfly of algebra}
 \xymatrix{
\frkg_1\ar[dd]\ar[dr]^{\kappa}&&\frkh_1\ar[dl]_{i_2}\ar[dd]\\
&\frkg_0\oplus\frkh_1\ar[dl]^{\sigma}\ar[dr]_{\rho}&\\
 \frkg_0&& \frkh_0. }
\end{equation} This correspondence of butterflies and morphisms is one-to-one, and  the morphism corresponding to a butterfly
is an equivalence of the corresponding Lie 2-algebras associated to the crossed modules if and only if the NW-SE sequence is
a short exact sequence (see \cite[Remark 5.5]{noohi:morphism}).
\begin{rmk}\label{rk:integration}
Butterfly provides a convenient tool to integrate morphisms.  Notice that a short exact sequence of Lie
algebras integrates to the short exact sequence of the corresponding
simply connected Lie groups. We can integrate the butterfly \eqref{eqn:butterfly of algebra} to a butterfly of crossed modules
of Lie groups by integrating all the Lie algebras to their simply
connected Lie groups. It is also easy to see that an equivalence of crossed module of Lie algebras integrates to an
equivalence of the corresponding crossed module of Lie groups. See Proposition 3.4 in
\cite{noohi:morphism} for more details.
\end{rmk}

\section{Integrating semidirect product Lie 2-algebras}\label{sec:integrating-lie2}

\subsection{Integrating the crossed module of Lie algebras
$\End(\V)$}

The Lie algebra structure on $\frkk_0$ given by (\ref{k0}) is very
clear. The Lie group $K_0$ defined by
\begin{equation}\label{group K0}
K_0=\{\Big(\begin{array}{cc}B_0&0\\0&B_1\end{array}\Big):\quad
B_0\in GL(V_0), ~B_1\in GL(V_1)\mbox{ such that }B_0\circ
\dM=\dM\circ B_1\}
\end{equation}
is a Lie group whose Lie algebra is $\frkk_0$.

The Lie algebra structure on $\frkk_1$ given by (\ref{k1}), i.e. the
Lie bracket $[\cdot,\cdot]_{\frkk_1}$ is less clear. We consider the
integration of the Lie algebra $\frkk_1$ and the Lie algebra
morphism $\delta$ by embedding $\frkk_1$ into the Lie algebra
$\gl(V_0\oplus V_1)$.

\begin{pro}\label{thm:int cm}
The set $K_1$  given by
\begin{equation}\label{group k1}
K_1=\{M\in\End(V_0,V_1)\quad\mbox{such that}\quad\Big(
\begin{array}{cc}I+\dM\circ M&0\\0&M\circ\dM+I\end{array}\Big)\in GL(V_0\oplus
V_1)\},
\end{equation}
is a Lie group with  multiplication \be\label{m in H1} M_1\cdot
M_2=M_1+ M_2 +M_1\circ\dM\circ M_2,\ee identity $0$ and inverse
$$
M^{-1}=-(I+M\circ \dM)^{-1} \circ M=-M\circ (I+ \dM\circ M)^{-1}.
$$
Its  Lie algebra is $\frkk_1$ given by \eqref{k1}. The Lie group
morphism $\int\delta$ from $K_1$ to $K_0$, given by \be\label{int
delta}
(\int\delta)(M)=\Big(\begin{array}{cc} I+\dM\circ M&0\\
0&M\circ\dM+I\end{array}\Big), \ee differentiates to the Lie algebra
morphism $\delta: \frkk_1\to \frkk_0$.  The action $\Phi$ of the Lie
group $K_0$ on $K_1$ given by
\begin{equation}\label{action of H0}
\Phi_{\Big(\begin{array}{cc}B_0&0\\0&B_1\end{array}\Big)}M=B_1\circ
M\circ B_0^{-1},\quad
\Big(\begin{array}{cc}B_0&0\\0&B_1\end{array}\Big)\in K_0, ~M\in
K_1,
\end{equation}
differentiates to the action $\phi$ of Lie algebra $\frkk_0$ on
$\frkk_1$ in \eqref{eq:phi}. Furthermore,
\begin{equation}\label{crossed module autv}
\Aut(\V)\triangleq(K_1,K_0,\int\delta,\Phi)
\end{equation}
is a crossed module of Lie groups whose differentiation is the
crossed module of Lie algebras $\End(\V)$ given by \eqref{crossed
module endv}.
\end{pro}
\pf It is not hard to see that the group structure of $K_1$ has the
property that the injective map $M \mapsto
\Big(\begin{array}{cc}I+\dM\circ M&0\\0&M\circ\dM+I\end{array}\Big)$
is a group morphism from $K_1$ to $GL(V_0\oplus V_1)$. It is obvious
that $0$ is the identity of the multiplication $(\ref{m in H1})$. It
is also not hard to see that
\begin{eqnarray*}
M\cdot(-(I+M\circ \dM)^{-1} \circ M)&=&M+(I+M\circ\dM)(-(I+M\circ
\dM)^{-1} \circ M)=0,\\
(-M\circ (I+ \dM\circ M)^{-1})\cdot M&=&(-M\circ (I+ \dM\circ
M)^{-1})(I+\dM\circ M)+M=0.
\end{eqnarray*}
Furthermore, the equality
$$
M+M\circ\dM\circ M=M\circ(I+\dM\circ M)=(I+M\circ\dM)\circ M
$$
yields that
$$
-(I+M\circ \dM)^{-1} \circ M=-M\circ (I+ \dM\circ M)^{-1}.
$$
We can view $K_1$ as a subgroup of $GL(V_0\oplus V_1)$.

Now we identify $\frkk_1$ and $K_1$ to their images in the
corresponding bigger matrix spaces. By exponentiating, we have
$$
\exp{\Big(\begin{array}{cc} \dM\circ A&0\\
&A\circ\dM\end{array}\Big)}=\Big(\begin{array}{cc} I+\dM\circ M&0\\
&M\circ\dM+I\end{array}\Big),
$$
where $$M=A+\frac{A\circ\dM\circ A}{2!}+\frac{A\circ\dM\circ
A\circ\dM\circ A}{3!}+\cdots.$$

Since the exponential map is an isomorphism near 0, by
\cite[Prop.3.18]{warner}, $K_1$ is a Lie sub-group of $GL(V_0\oplus
V_1)$ whose Lie algebra is $\frkk_1$.   It is not hard to see that
the differentiation of $\int \delta$ given by (\ref{int delta}) is
the $\delta$ given by (\ref{eqn:differential}). Since the Lie
algebra action is given by the commutator, it follows that the group
action is given by the adjoint action, which turns out to be
(\ref{action of H0}). In other words, \eqref{action of H0} is
constructed so that its differentiation is the $\phi$ in
\eqref{eq:phi}.

Finally, we prove that $\Aut(\V)=(K_1,K_0,\int\delta,\Phi)$ is a
crossed module of Lie groups. First of all, we have
\begin{eqnarray*}
(\int\delta)\Phi_{\Big(\begin{array}{cc}B_0&0\\0&B_1\end{array}\Big)}M&=&(\int\delta)(B_1\circ
M\circ B_0^{-1})\\
&=&\Big(\begin{array}{cc}I+\dM\circ B_1\circ M\circ
B_0^{-1}&0\\0&I+B_1\circ M\circ B_0^{-1}\circ\dM\end{array}\Big)\\
&=&\Big(\begin{array}{cc}B_0&0\\0&B_1\end{array}\Big)
(\int\delta)(M)\Big(\begin{array}{cc}B_0&0\\0&B_1\end{array}\Big)^{-1}.
\end{eqnarray*}
Furthermore, we have
\begin{eqnarray*}
\Phi_{(\int\delta)(M)}(M^\prime)&=&\Phi_{\Big(\begin{array}{cc}I+\dM\circ
M&0\\0&I+ M\circ\dM\end{array}\Big)}M^\prime\\
&=&(I+ M\circ\dM)\circ M^\prime \circ(I+\dM\circ M)^{-1}.
\end{eqnarray*}
On the other hand, using the facts that $ M^{-1}= -M\circ(I+\dM\circ
M)^{-1}$ and $I=(I+\dM\circ M)\circ(I+\dM\circ M)^{-1}$, we have
\begin{eqnarray*}
&&M\cdot M^\prime\cdot M^{-1}\\&=&(M+M^\prime+M\circ\dM\circ
M^\prime)\cdot(-M\circ(I+\dM\circ M)^{-1})\\
&=&M+M^\prime+M\circ\dM\circ M^\prime-M\circ(I+\dM\circ
M)^{-1}-M\circ\dM\circ M\circ(I+\dM\circ M)^{-1}\\
&&-M^\prime\circ\dM \circ M\circ(I+\dM\circ M)^{-1}-M\circ\dM\circ
M^\prime\circ \dM\circ M\circ(I+\dM\circ M)^{-1}\\
&=&M^\prime\circ(I-\dM\circ M\circ(I+\dM\circ
M)^{-1})+M\circ\dM\circ
M^\prime\circ(I-\dM\circ M\circ(I+\dM\circ M)^{-1})\\
&=&M^\prime\circ(I+\dM\circ M)^{-1}+M\circ\dM\circ M^\prime\circ
(I+\dM\circ M)^{-1}.
\end{eqnarray*}
Thus we have
$$
\Phi_{(\int\delta)(M)}(M^\prime)=M\cdot M^\prime\cdot M^{-1},
$$
which implies that $\int \delta$ satisfies the Peiffer identity.
Therefore,  $\Aut(\V)=(K_1,K_0,\int\delta,\Phi)$ is a crossed module
of Lie groups. \qed

\subsection{Semidirect product Lie 2-groups--strict case}
Let $\HH=(H_1,H_0,t,\Phi)$ and $\G=(G_1,G_0,t,\Phi)$ be crossed
modules of Lie groups, and let
$\V=(V_1\stackrel{\dM}{\longrightarrow}V_0)$ be  a 2-term complex of
vector spaces. In the following, $x,y,z$ are elements of $H_0$,
$a,b,c$ are elements of $H_1$,  $\xi,\eta,\gamma$ are elements of
$ V_0$, and $m,n,p,q$ are elements of $V_1$. We omit the group multiplication $\cdot$ in $H_i$'s, $K_i$'s, and $G_i$'s, and denote the composition of morphisms by $\circ$.
\begin{defi}\label{def:m of crossed modules of groups}
A strict morphism from $\HH$ to $\G$ is a pair $(F_1,F_0)$, where
$F_i:H_i\longrightarrow G_i$ are morphisms of Lie groups, such that
$$
F_0\circ t=t\circ F_1,\quad F_1(\Phi_xa)=\Phi_{F_0(x)}F_1(a).
$$
\end{defi}
 Let $(\Psi_1,\Psi_0)$ be a strict morphism  of crossed modules of Lie
groups from $\HH$ to $\Aut(\V)$.  By Definition \ref{def:m of
crossed modules of groups} and Proposition \ref{thm:int cm}, we have
\begin{eqnarray}\label{lem:adj}
\Psi_1(\Phi_xa)=\Psi_0(x)\circ\Psi_1(a)\circ(\Psi_0(x))^{-1},
\end{eqnarray}
and
\begin{eqnarray*}
(\int\delta)\Psi_1(a)=\Big(\begin{array}{cc} I+\dM\circ \Psi_1(a)&0\\
&\Psi_1(a)\circ\dM+I\end{array}\Big)=\Psi_0(t(a)).
\end{eqnarray*}
More precisely, for any $\xi\in V_0$ and $m\in V_1$, we have
\begin{eqnarray}
\label{homotopy 1}\xi+\dM\circ \Psi_1(a)(\xi)&=&\Psi_0(t(a))(\xi),\\
\label{homotopy 2}m+ \Psi_1(a)\circ\dM m&=&\Psi_0(t(a))(m).
\end{eqnarray}

\begin{lem}\label{lem:m}
For any $ a,b\in H_1$, we have
$$
\Psi_1(ab)=\Psi_1(a)+\Psi_0(t(a))\circ\Psi_1(b).
$$
\end{lem}
\pf By (\ref{m in H1}) and the fact that $\Psi_1$ is a morphism, we
have
\begin{eqnarray*}
\Psi_1(ab)&=&\Psi_1(a)\Psi_1(b)=\Psi_1(a)+\Psi_1(b)+\Psi_1(a)\circ\dM\circ\Psi_1(b)\\
&=&\Psi_1(a)+(I+\Psi_1(a)\circ\dM)\circ\Psi_1(b)\\
&=&\Psi_1(a)+(\int\delta)\Psi_1(a)\circ\Psi_1(b)\\
&=&\Psi_1(a)+\Psi_0(t(a))\circ\Psi_1(b).  \qed
\end{eqnarray*}

\begin{lem}\label{lem:inverse}
For any $a \in H_1$, we have
\begin{equation}\label{inverse}
\Psi_0(t(a^{-1}))\circ\Psi_1(a)+\Psi_1(a^{-1})=0.
\end{equation}
\end{lem}
\pf By (\ref{m in H1}), we have
$$
0=\Psi_1(a^{-1}a)=\Psi_1(a^{-1})\Psi_1(a)=\Psi_1(a)+\Psi_1(a^{-1})+\Psi_1(a^{-1})\circ\dM\circ\Psi_1(a).
$$
Thus,\begin{eqnarray*}
\Psi_0(t(a^{-1}))\circ\Psi_1(a)+\Psi_1(a^{-1})&=&(\int\delta)(\Psi_1(a^{-1}))\circ\Psi_1(a)-\Psi_1(a)-\Psi_1(a^{-1})\circ\dM\circ\Psi_1(a)\\
&=&(I+\Psi_1(a^{-1})\circ\dM)\circ\Psi_1(a)-\Psi_1(a)-\Psi_1(a^{-1})\circ\dM\circ\Psi_1(a)\\
&=&0.\qed
\end{eqnarray*}

\begin{thm}\label{thm:semidirect product of strict}
Given a strict morphism $(\Psi_1,\Psi_0)$ of crossed modules of Lie
groups from $(H_1,H_0,t,\Phi)$ to $\Aut(\V)$ given by \eqref{crossed
module autv}, there is a strict Lie 2-group\begin{equation}
\label{eq:smdr-prd-gp} \begin{array}{c}
H_0\times H_1\times V_0\times V_1\\
\vcenter{\rlap{s }}~\Big\downarrow\Big\downarrow\vcenter{\rlap{t }}\\
H_0\times V_0,
 \end{array}\end{equation}
 in which the source and the target maps are given by
 \begin{eqnarray*}
s(x,a,\xi,m)&=&(x,\xi),\\
t(x,a,\xi,m)&=&(t(a)x,\xi+\dM m),
 \end{eqnarray*}
the vertical multiplication $\cdot_\ve$ is given by
$$
(y,b,\eta,n)\cdot_\ve(x,a,\xi,m) =(x,ba,\xi,m+n),\quad
y=t(a)x,~\eta=\xi+\dM m,
$$
the horizontal multiplication $\cdot_\mathrm{h}$ of morphisms is
defined by
\begin{eqnarray*} &(x,a,\xi,m)\cdot_\mathrm{h}(y,b,\eta,n)=\Big(xy,a\Phi_{x}b,\xi+
\Psi_0(x)(\eta),m+\Psi_0(t(a)x)(n)+\Psi_1(a)\circ
\Psi_0(x)(\eta)\Big),\end{eqnarray*} and the horizontal multiplication
$\cdot_\mathrm{h}$ of objects is defined by
  $$
(x,\xi)\cdot_\mathrm{h}(y,\eta)=(xy,\xi+\Psi_0(x)(\eta)).
  $$
The identities of arrows
and objects are $(1_{H_0},1_{H_1},0,0)$ and $(1_{H_0},0)$  respectively. The inverse of $(x,a,\xi,m)$ with respect to $\cdot_h$ is
$$\big(x^{-1},\Phi_{x^{-1}}a^{-1},-\Psi_0(x^{-1})(\xi),-\Psi_0((t(a)x)^{-1})(m)+\Psi_0((t(a)x)^{-1})\circ\Psi_1(a)(\xi)\big).$$
\end{thm}
We call $(\Psi_1, \Psi_0)$ a {\em representation} on $\V$ of the
crossed module $(H_1,H_0,t,\Phi)$, and  the strict Lie 2-group
\eqref{eq:smdr-prd-gp} the {\em semidirect product} of
$(H_1,H_0,t,\Phi)$ with this representation. \vspace{3mm}

\pf Obviously, the horizontal multiplication respects the source
map. By (\ref{homotopy 1}) and (\ref{homotopy 2}), the horizontal
multiplication also respects  the target map.

To see the horizontal multiplication is indeed a functor, we need to
show that
$$
\Big((t(a)x,b,\xi+\dM m,n)\cdot_\ve(x,a,\xi,m)\Big)\cdot_\mathrm{h}
\Big((t(c)y,d,\eta+\dM p,q)\cdot_\ve(y, c,\eta,p)\Big)
$$
is equal to $$ \Big((t(a)x,b,\xi+\dM
m,n)\cdot_\mathrm{h}(t(c)y,d,\eta+\dM p,q)\Big)\cdot_\ve
\Big((x,a,\xi,m)\cdot_\mathrm{h}(y,c,\eta,p)\Big),
$$
where $x,y\in H_0,~ a,b,c,d\in H_1, ~\xi,\eta,\gamma\in V_0$, and
$m,n,p,q\in V_1$. By straightforward computations, the first
expression is equal to
\begin{equation}\label{temp1}
\Big(xy,ba\Phi_x(dc),\xi+\Psi_0(x)(\eta),m+n+\Psi_0(t(ba)x)(p+q)+\Psi_1(ba)\circ\Psi_0(x)(\eta)\Big).
\end{equation}
The second expression is equal to
\begin{eqnarray*}\label{temp2}
&&\Big(xy,(b\Phi_{t(a)x}d)(a\Phi_xc),\xi+\Psi_0(x)(\eta),\\
&&m+n+\Psi_0(t(a)x)(p)+\Psi_0(t(ba)x)(q)+\Psi_1(a)\circ\Psi_0(x)(\eta)+\Psi_1(b)\circ\Psi_0(t(a)x)(\eta+\dM
p)\Big).
\end{eqnarray*}
By the definition of crossed modules, it is not hard to see that
$$
(b\Phi_{t(a)x}d)(a\Phi_xc)=b\Phi_{t(a)}(\Phi_xd)a\Phi_xc=ba(\Phi_xd)a^{-1}a\Phi_xc=ba\Phi_x(dc).
$$
Thus we only need to prove that
\begin{eqnarray}
\nonumber&&\Psi_0(t(ba)x)(p)+\Psi_1(ba)\circ\Psi_0(x)(\eta)\\
&=&\Psi_0(t(a)x)(p)+\Psi_1(a)\circ\Psi_0(x)(\eta)+\Psi_1(b)\circ\Psi_0(t(a)x)(\eta+\dM
p).
\end{eqnarray}
By (\ref{m in H1}), we have
\begin{equation}\label{temp3} \Psi_1(ba)=\Psi_1(b)+\Psi_1(a)+\Psi_1(b)\circ\dM
\circ\Psi_1(a).
\end{equation}
So it suffices to prove that \begin{eqnarray*}
\Psi_0(t(ba)x)(p)&=&\Psi_0(t(a)x)(p)+\Psi_1(b)\circ\Psi_0(t(a)x)(\dM p),\\
\text{and} \quad \Psi_1(b)\circ\Psi_0(t(a)x)(\eta)&=&\Psi_1(b)\circ\Psi_0(x)(\eta)+\Psi_1(b)\circ\dM
\circ\Psi_1(a)\circ\Psi_0(x)(\eta),
\end{eqnarray*}
which hold if and only if for any $\xi\in V_0, m\in V_1$, we have
\begin{eqnarray*}
(\int\delta)(\Psi_1(b))(m)&=&I+\Psi_1(b)\circ\dM m,\\
(\int\delta)(\Psi_1(a))(\xi)&=&I+\dM\circ\Psi_1(a)(\xi).
\end{eqnarray*}
These are exactly (\ref{homotopy 1}) and (\ref{homotopy 2}). Thus the
horizontal multiplication $\cdot_\mathrm{h}$ is indeed  a functor.

To see that the horizontal multiplication $\cdot_{\mathrm h}$ is
strictly associative, we only need to check it on the
  space of arrows, i.e. we need to verify that
\begin{equation}
\big((x,a,\xi,m)\cdot_\hh(y,b,\eta,n)\big)\cdot_\hh(z,c,\gamma,p)=(x,a,\xi,m)\cdot_\hh\big((y,b,\eta,n)\cdot_\hh(z,c,\gamma,p)\big).
\end{equation}
By straightforward computations, we obtain that the left hand side
is equal to
\begin{eqnarray*}
&\Big(xyz,a(\Phi_xb)(\Phi_{xy}c),\xi+\Psi_0(x)(\eta)+\Psi_0(xy)(\gamma),\\&m+\Psi_0(t(a)x)(n)+\Psi_1(a)\circ\Psi_0(x)(\eta)+\Psi_0(t(a\Phi_xb)xy)(p)
+\Psi_1(a\Phi_xb)\circ\Psi_0(xy)(\gamma)\Big),
\end{eqnarray*}
and the right hand side is equal to
\begin{eqnarray*}
&\Big(xyz,a\Phi_x(b\Phi_yc),\xi+\Psi_0(x)(\eta)+\Psi_0(x)\circ\Psi_0(y)(\gamma),m+\Psi_0(t(a)x)(n)+
\Psi_1(a)\circ\Psi_0(x)(\eta)\\&+\Psi_0(t(a)x)\circ\Psi_0(t(b)y)(p)
+\Psi_0(t(a)x)\circ\Psi_1(b)\circ\Psi_0(y)(\gamma)+\Psi_1(a)\circ\Psi_0(x)\circ\Psi_0(y)(\gamma)\Big).
\end{eqnarray*}
 Since $\Psi_0$ is a morphism of Lie groups and $\Phi$ acts as an automorphism, we only need to show that
 $$
\Psi_1(a\Phi_xb)\circ\Psi_0(xy)(\gamma)=\Psi_0(t(a)x)\circ\Psi_1(b)\circ\Psi_0(y)(\gamma)+\Psi_1(a)\circ\Psi_0(x)\circ\Psi_0(y)(\gamma).
 $$
 Since
 $\Psi_1(a\Phi_xb)=\Psi_1(a)+\Psi_1(\Phi_xb)+\Psi_1(a)\circ\dM\circ\Psi_1(\Phi_xb)$,
 it is equivalent to
\begin{eqnarray*}
&&\Psi_1(\Phi_xb)\circ\Psi_0(x)\circ\Psi_0(y)(\gamma)+\Psi_1(a)\circ\dM\circ\Psi_1(\Phi_xb)\circ\Psi_0(x)\circ\Psi_0(y)(\gamma)\\
&&=\Psi_0(t(a)x)\circ\Psi_1(b)\circ\Psi_0(y)(\gamma),
\end{eqnarray*}
which holds by (\ref{homotopy 2}).

Finally, we show that $(1_{H_0},1_{H_1},0,0)$ and $(1_{H_0},0)$ are
identities of arrows and objects respectively. It is straightforward
to see that for any $(x,a,\xi,m)$, we have
\begin{eqnarray*}
&&(x,a,\xi,m)\big(x^{-1},\Phi_{x^{-1}}a^{-1},-\Psi_0(x^{-1})(\xi),-\Psi_0((t(a)x)^{-1})(m)+\Psi_0((t(a)x)^{-1})\circ\Psi_1(a)(\xi)\big)\\
&&=(1_{H_0},1_{H_1},0,0).
\end{eqnarray*}
On the other hand, we have
\begin{eqnarray*}
&&\big(x^{-1},\Phi_{x^{-1}}a^{-1},-\Psi_0(x^{-1})(\xi),-\Psi_0((t(a)x)^{-1})(m)+\Psi_0((t(a)x)^{-1})\circ\Psi_1(a)(\xi)\big)(x,a,\xi,m)\\
&&=\big(1_{H_0},1_{H_1},0,\Psi_0((t(a)x)^{-1})\circ\Psi_1(a)(\xi)+\Psi_1(\Phi_{x^{-1}}a^{-1})\circ\Psi_0(x^{-1})(\xi)\big).
\end{eqnarray*}
To see that
$\big(x^{-1},\Phi_{x^{-1}}a^{-1},-\Psi_0(x^{-1})(\xi),-\Psi_0((t(a)x)^{-1})(m)+\Psi_0((t(a)x)^{-1})\circ\Psi_1(a)(\xi)\big)$
is the inverse of $(x,a,\xi,m)$, we need to prove that
\begin{equation}
\Psi_0(t(a^{-1}))\circ\Psi_1(a)(\xi)+\Psi_1(a^{-1})(\xi)=0.
\end{equation}
This is exactly Lemma \ref{lem:inverse}. The proof is finished.
\qed\vspace{3mm}


\begin{rmk} \label{cor:cm group} The strict Lie 2-group constructed in the above theorem
  can be viewed as the semidirect product $\HH \ltimes \V$. We
construct the formula from its infinitesimal counter part $\h
\ltimes \V$ in Example \ref{ep:semidirect}. In fact, our Lie 2-group
corresponds to the crossed module of Lie groups: $\big(H_1\ltimes
V_1,H_0\ltimes V_0,t\times\dM,\Phi\big)$, where
$t\times\dM:H_1\ltimes V_1\longrightarrow H_0\ltimes V_0$ is given
by
$$
(t\times\dM)(a,m)=(t(a),\dM m),
$$
the group structure on $H_1\ltimes V_1$ is given by
$$
(a,m)(b,n)=\big(ab,m+\Psi_0(t(a))(n)\big),
$$
the group structure on $H_0\ltimes V_0$ is the semidirect product,
i.e.
$$
(x,\xi)(y,\eta)=(xy,\xi+\Psi_0(x)(\eta)),
$$
and the action of $H_0\ltimes V_0$ on $H_1\ltimes V_1$ is given by
\begin{equation}\label{action of phi}
\Phi_{(x,\xi)}(a,m)=\big(\Phi_xa,\Psi_0(x)(m)-\Psi_0(x)\circ\Psi_1(a)\circ\Psi_0(x^{-1})(\xi)\big).
\end{equation}

\end{rmk}

\begin{rmk}\label{rk:simply connected}
Lie's II and III theorems hold for crossed modules, thus strict Lie
2-algebras. Given a crossed module of Lie algebras $\h$,
theoretically, there is a unique crossed module of Lie groups
$\HH=(H_1, H_0, t, \Phi)$ such that its differentiation is $\h$ and
$H_i$'s are connected and simply connected. We call  $\HH$ the {\bf
simply connected integration} of $\h$. In Theorem
\ref{thm:semidirect product of strict} and Remark \ref{cor:cm
group}, we explicitly constructed the integration of a semidirect
product Lie 2-algebra under a strict morphism.

Moreover, given any
crossed module of Lie groups $\G$ whose differentiation is $\g$,  a
morphism of crossed module of Lie algebras
$(\psi_1,\psi_0):\h\longrightarrow \g$ can be integrated to a
morphism of crossed module of Lie groups
$(\Psi_1,\Psi_0):\HH\longrightarrow\G$ where $\HH$ is the simply
connected integration of $\h$. See \cite{noohi:morphism} for more
details about the integration of morphisms of crossed modules.
\end{rmk}

Finally, we summarize what we obtain:

\begin{cor}\label{thm:semidirect product of strict 1}
Given a strict morphism $(\psi_1,\psi_0)$ from $\h$ to $\End(\V)$,
let $\HH$ be the simply connected integration of $\h$ and
$(\Psi_1,\Psi_0): \HH \to Aut(\V)$ be the integration (see Remark
\ref{rk:simply connected}) of $(\psi_1,\psi_0)$. Then the strict Lie
2-group given in Theorem \ref{thm:semidirect product of
  strict} is the
simply connected integration of $\h\ltimes \V$.
\end{cor}

If in Corollary \ref{cor:cm group}, we take $\Psi_0,\Psi_1$ to be the
identity map, then we obtain the crossed module of Lie groups
$$(K_1\ltimes V_1,K_0\ltimes
V_0,\int\delta\times\dM,\Phi),$$ which plays the role of
$GL(V)\ltimes V$ in the classical case of a vector space $V$ acted
upon by $GL(V)$, where $K_1$ and $K_0$ are given by (\ref{group k1})
and (\ref{group K0}). If in Corollary \ref{thm:semidirect product of
strict 1}, we take $\psi_0,\psi_1$ to be the identity map, then we
also obtain a crossed module of Lie algebras
$$
(\frkk_1\ltimes V_1,\frkk_0\ltimes V_0,\delta\times\dM,\phi),
$$
which plays the role of $\gl(V)\ltimes V$ in the classical case of a vector space $V$ acted upon by $\gl(V)$, where $\frkk_1$ and $\frkk_0$ are given by
(\ref{k1}) and (\ref{k0}).

\subsection{Non-strict case}

\begin{defi}
A strict Lie 2-group is called an integration of  a Lie 2-algebra if
its differentiation is a Lie 2-algebra equivalent to the given Lie
2-algebra.
\end{defi}

\begin{rmk}
When one studies the  integration of the string Lie
2-algebras \cite{baez:str-gp, henriques,
schommer:string-finite-dim}, the models one finds are (only equivalent) not the same, as 2-groups. Equivalence of
2-groups corresponds to equivalence of Lie 2-algebras on the
infinitesimal level. Thus, motivated by the examples arising from the integration of string Lie 2-algebras, we define our integration up to equivalence
as above.
\end{rmk}

The final aim of this paper is to integrate the Lie 2-algebra
$\frkg\ltimes\V$, which is the semidirect product of a Lie algebra
$\frkg$ with a 2-term representation up to homotopy $\V$. Recall
that a 2-term representation up to homotopy is a non-strict morphism
$(\mu,\nu)$ from $\frkg$ to $\End(\V)$. Moreover, in the last
subsection we study the integration of the Lie 2-algebra which is
the semidirect product of $\h$ with $\V$ via a strict morphism. Thus
we give the integration of the Lie 2-algebra $\frkg\ltimes\V$ by
strictifying the non-strict morphism $(\mu,\nu)$ via the
corresponding butterfly.

For any butterfly $\frke$ from  $\g$ to $\h$,
\begin{equation}\label{diagrambutterfly0}
 \xymatrix{
\frkg_1\ar[dd]\ar[dr]^{\kappa}&&\frkh_1\ar[dl]_{\iota}\ar[dd]\\
&\frke\ar[dl]^{\sigma}\ar[dr]_{\rho}&\\
 \frkg_0&& \frkh_0, }
\end{equation}
we obtain a crossed module of Lie algebras
\begin{equation}\label{eqn:crossed
module}(\frkg_1\times_{\frkg_0}\frke,\frke,dt,\phi),
\end{equation}where $\frkg_1\times_{\frkg_0}\frke$ is short for the fibre product $\frkg_1\times_{dt, \frkg_0, \sigma} \frke$.
The Lie algebra structure on $\frkg_1\times_{\frkg_0}\frke$ is given
by
$$
[(A_1,e_1),(A_2,e_2)]=([A_1,A_2],[e_1,e_2]),
$$
the Lie algebra morphism
$dt:\frkg_1\times_{\frkg_0}\frke\longrightarrow\frke$ is defined by
$$
dt(A,e)=e,\quad\forall ~A\in\frkg_1,~e\in \frke,~dt(A)=\sigma(e),
$$
and the action $\phi$ of $\frke$ on $\frkg_1\times_{\frkg_0}\frke$
is given by
$$
\phi_e(A,e_1)=(\phi_{\sigma(e)}A,[e,e_1]).
$$

Furthermore, there are two strict morphisms of crossed modules of
Lie algebras, which are $(\pr_1,\sigma)$ and $(\psi_1,\psi_0)$,
\begin{equation}\label{diagrambutterfly1}
 \xymatrix{&\frkg_1\times_{\frkg_0}\frke\ar[dd]^{dt}\ar[dl]^{\pr_1}\ar[dr]_{\psi_1}&\\
\frkg_1\ar[dd]\ar[dr]^{\kappa}&&\frkh_1\ar[dl]_{\iota}\ar[dd]\\
&\frke\ar[dl]^{\sigma}\ar[dr]_{\rho=\psi_0}&\\
 \frkg_0&& \frkh_0, }
\end{equation}
where $\psi_0=\rho:\frke\longrightarrow \frkh_0$ and
$\psi_1:\frkg_1\times_{\frkg_0}\frke\longrightarrow\frkh_1$ is
defined by
\begin{equation}\label{eqn:map}
\psi_1(A,e)=e-\kappa(A),\quad\forall~A\in\frkg_1,~e\in \frke.
\end{equation}

The following conclusion is straightforward.
\begin{pro}Given a butterfly \eqref{diagrambutterfly0}, the map $(\pr_1,\sigma)$ constructed in the diagram \eqref{diagrambutterfly1}
is a strict morphism from the crossed module
$(\frkg_1\times_{\frkg_0}\frke,\frke,dt,\phi)$  to $\g$. The map
$(\pr_1,\sigma)$ is also an equivalence of the underlying 2-term
complexes. Thus, the corresponding Lie 2-algebras are equivalent.
\end{pro}

At the end of this section, we focus on non-strict  morphisms
$(\mu,\nu)$ from a Lie algebra $\frkg$, which is viewed as a trivial
crossed module of Lie algebras, to the crossed module of Lie
algebras $\End(\V)=(\frkk_1, \frkk_0, dt, \phi)$, which is given by
(\ref{crossed module endv}). Obviously, the corresponding butterfly
is $\frke=\frkg\oplus\frkk_1$ given by (\ref{eqn:butterfly of
algebra}). The crossed module of Lie algebras given by
(\ref{eqn:crossed module}) turns out to be
$(\sigma^{-1}(0),\frkg\oplus \frkk_1,dt,\phi)$. Moreover, it is
straightforward to see that $\sigma^{-1}(0)$ is exactly $\frkk_1$.
Thus the corresponding crossed module of Lie algebras is
\begin{equation}\label{eq:g-pull-back}
(\frkk_1,\frkg\oplus \frkk_1,i_2,\ad),
\end{equation}
where the Lie algebra structure on $\frkg\oplus\frkk_1$ is given by
$$
[(X,A),(Y,B)]=\big([X,Y],[\mu(X),B]+[A,\mu(Y)]+[A,B]_{\frkk_1}-\nu(X,Y)\big),
$$
and the adjoint action $\ad$ is given by
$$
\ad_{(X,A)}B=\phi_{\mu(X)}B+[A,B]_{\frkk_1}=[\mu(X),B]+A\circ\dM\circ
B-B\circ\dM\circ A.
$$
Furthermore, the Lie algebra morphism $\psi_1:\frkk_1\longrightarrow
\frkk_1$ given by (\ref{eqn:map}) is exactly the identity map $\Id$
and the Lie algebra morphism
$\psi_0=\rho:\frkg\oplus\frkk_1\longrightarrow\frkk_0$ is given by
$$
\psi_0(X,A)=\mu(X)+\delta(A).
$$Thus $(\mu, \nu)$ strictifies to a strict morphism
$(\psi_0, \psi_1)$.  The
semidirect product of \eqref{eq:g-pull-back} with $(\psi_1, \psi_0)$
as in Corollary \ref{thm:semidirect product of strict 1} is
\begin{equation}\label{eq:sm-k}(\frkk_1\ltimes
V_1,(\frkg\oplus\frkk_1)\ltimes V_0,i_2\times\dM,\widetilde{\phi}),
\end{equation}
 where the
action $\widetilde{\phi}$ is given by
$$
\widetilde{\phi}_{(X,A,\xi)}(B,m)=\big(\ad_{(X,A)}B,(\mu(X)+\delta(A))(m)-B\xi\big).
$$

\begin{pro}\label{thm:equivalent}
With the above notations, given a (non-strict)  morphism of crossed
modules $(\mu, \nu): \frkg \to \End(\V)$, we have a strict morphism
$(\psi_1,\psi_0)=(\Id,\rho)$  from the crossed module defined by
\eqref{eq:g-pull-back} to $\End(\V)$. The semidirect product Lie
2-algebra \eqref{eq:sm-k} is equivalent to the semidirect product
Lie 2-algebra $\frkg\ltimes\V$.
\end{pro}
\pf  We only need to show that the Lie 2-algebras $(\frkk_1\ltimes
V_1,(\frkg\oplus\frkk_1)\ltimes V_0,i_2\times\dM,\widetilde{\phi})$
and $\frkg\ltimes \V$ are equivalent. Define
$f_0:(\frkg\oplus\frkk_1)\ltimes V_0\longrightarrow \frkg\ltimes
V_0$ by
$$
f_0(X,A,\xi)=(X,\xi),\quad\forall~(X,A)\in\frkg\oplus\frkk_1,~\xi\in
V_0,
$$
define $f_1:\frkk_1\ltimes V_1\longrightarrow V_1$ by
$$
f_1(A,m)=m,\quad\forall~A\in\frkk_1,~m\in V_1.
$$
Obviously, $(f_0,f_1)$  respects the differential, i.e
$$
(0+\dM)\circ f_1=f_0\circ (i_2\times\dM).
$$
Define $f_2:\wedge^2((\frkg\oplus\frkk_1)\ltimes V_0)\longrightarrow
V_1$ by
$$
f_2((X,A,\xi),(Y,B,\eta))=A\eta-B\xi.
$$
By straightforward computations, we have \begin{eqnarray*}
&&f_0[(X,A,\xi),(Y,B,\eta)]-[f_0(X,A,\xi),f_0(Y,B,\eta)]\\&=&
\rho(X,A)(\eta)-\rho(Y,B)(\xi)-\big(\mu(X)(\eta)-\mu(Y)(\xi)\big)\\
&=&\dM \circ A\eta-\dM \circ B\xi\\&=&\dM f_2((X,A,\xi),(Y,B,\eta)).
 \end{eqnarray*}
On the other hand,
\begin{eqnarray*}
f_1(\phi_{(X,A,\xi)}(B,m))&=&f_1\big(\ad_{(X,A)}B,\rho(X,A)(m)-\psi_1(B)(\xi)\big)\\&=&(\mu(X)+\delta(A))(m)-B\xi\\
&=&\mu(X)(m)+A\circ\dM m-B\xi,\\
 \phi_{f_0(X,A,\xi)}f_1(B,m)&=&\phi_{(X,\xi)}m=\mu(X)(m).
\end{eqnarray*}
Thus
\begin{eqnarray*}
f_1(\phi_{(X,A,\xi)}(B,m))-
 \phi_{f_0(X,A,\xi)}f_1(B,m)&=&A\circ\dM m-B\xi\\
 &=&f_2((X,A,\xi),(i_2\times\dM)(B,m)).
\end{eqnarray*}
Finally, it is straightforward to obtain that
\begin{eqnarray*}
&&f_2\big([(X,A,\xi),(Y,B,\eta)],(Z,C,\gamma)\big)+c.p.\\&=&[f_0(X,A,\xi),f_2((Y,B,\eta),(Z,C,\gamma))]+c.p.\\
&=&\mu(X)(B\gamma-C\eta)+\mu(Y)(C\xi-A\gamma)+\mu(Z)(A\eta-B\xi),
\end{eqnarray*}
which implies that $(f_0,f_1,f_2) $ is a Lie 2-algebra morphism.

Furthermore, it is obvious that $(f_0,f_1) $ also induces an
equivalence of the underlying complexes of vector spaces. Thus the
 Lie 2-algebra
$(\frkk_1\ltimes V_1,(\frkg\oplus\frkk_1)\ltimes
V_0,i_2\times\dM,\widetilde{\phi})$ given by \eqref{eq:sm-k} is
equivalent to the Lie 2-algebra $\frkg\ltimes\V$. \qed\vspace{2mm}

By Corollary \ref{thm:semidirect product of strict 1} and Theorem
\ref{thm:equivalent}, we have
\begin{cor} \label{thm:integration}
The strict Lie 2-group  corresponding to the simply connected
integration of the crossed module of Lie algebras \eqref{eq:sm-k} is
an integration of the Lie 2-algebra $\frkg\ltimes\V$.
\end{cor}

\begin{rmk}
In Corollary \ref{thm:integration}, we construct a strict Lie
2-group integrating the Lie 2-algebra $\frkg\ltimes\V$. In
\cite[Theorem 4.4]{shengzhu1}, the Lie 2-group constructed to
integrate the string type Lie algebra $\mathbb
R\longrightarrow\frkg\oplus \frkg^*$ is not strict since the
associator is not trivial. Thus this non-strict Lie 2-group is
equivalent to a strict Lie 2-group, which is the integration result
given in Corollary \ref{thm:integration}.
\end{rmk}

A sub-Lie algebra of a Lie algebra is defined by an injective Lie
algebra morphism. Similarly, we have:
\begin{defi}\label{def:subLie2}
Let $f:L\to W$ be a Lie 2-algebra morphism as in Definition
\ref{def:m of 2-alg}. Then $(L, f)$ is called a sub-Lie-2-algebra of
$W$ if $f_0$ and $f_1$ are injective. When the inclusion map $f$ is
obvious, we also call that $L$ is a sub-Lie-2-algebra of $W$.
\end{defi}
\begin{rmk} When $f_0$ and $f_1$ are injective, the linear category corresponding to $L$ is a subcategory of the one corresponding to $W$.
In \cite{baez:string}, the notion of Lie sub-2-algebra  is used and
it is also called a sub-Lie-2-algebra. We did not find its exact
definition. It seems that \cite{baez:string} requires a Lie
sub-2-algebra to be a sub-complex which is closed under the
differential and the brackets. If this is the case, a Lie
sub-2-algebra of a strict Lie 2-algebra must be strict. However, our
definition of sub-Lie 2-algebra is much weaker as we will see below.
\end{rmk}

The equivalence of $\frkg \ltimes \V$ and \eqref{eq:sm-k} can also
provided by  the canonical inclusion map $(i_2,i_1\times\Id)$,
\begin{eqnarray*}
i_2(m)&=&(0,m),\quad\forall~ m\in V_1,\\
(i_1\times\Id) (X,\xi)&=&(X,0,\xi),\quad\forall~ (X,\xi)\in
\frkg\oplus V_0.
\end{eqnarray*}
By straightforward computations, we have
\begin{eqnarray*}
[(i_1\times\Id)(X,\xi),( i_1\times\Id)(Y,\eta)]&=&\big([X,Y],-\nu(X,Y),\mu(X)(\eta)-\mu(Y)(\xi)\big),\\
(i_1\times\Id)[(X,\xi),(Y,\eta)]&=&\big([X,Y],0,\mu(X)(\eta)-\mu(Y)(\xi)\big).
\end{eqnarray*}
Therefore,
$$
(i_1\times\Id)[(X,\xi),(Y,\eta)]-[(
i_1\times\Id)(X,\xi),(i_1\times\Id)(Y,\eta)]=(i_2\times\dM)(\nu(X,Y),0).
$$
Define $\widetilde{\nu}:\wedge^2(\frkg\oplus V_0)\longrightarrow
\frkk_1\ltimes V_1$ by
$$
\widetilde{\nu}((X,\xi),(Y,\eta))=(\nu(X,Y),0).
$$
The following proposition is straightforward.
\begin{pro}\label{pro:sub-lie-2}
Given a nonstrict  morphism $(\mu,\nu)$ from a Lie algebra $\frkg$
to  $\End(\V)$, the map $(i_2,i_1\times \Id,\widetilde{\nu})$ is a
Lie 2-algebra morphism from $\frkg\ltimes \V$ to the strict Lie
2-algebra \eqref{eq:sm-k}. Consequently, the Lie 2-algebra
$\frkg\ltimes \V$ is a sub-Lie-2-algebra of the strict Lie 2-algebra
\eqref{eq:sm-k}.
\end{pro}

\section{Integrating the omni-Lie algebra $\gl(V)\oplus V$} \label{sec:omni}
The notion of omni-Lie algebra was introduced by A. Weinstein in
\cite{weinstein:omni} to characterize Lie algebra structures on a
vector space $V$. On the direct sum space $\gl(V)\oplus V$, the
nondegenerate symmetric $V$-valued pairing $\pair{\cdot,\cdot}$ is
given by
$$
\pair{(A,u),(B,v)}=\half(Av+Bu),
$$
and the bracket operation $\Courant{\cdot,\cdot}$ is given by
\begin{equation}\label{def:bracketomni}
\Courant{(A,u),(B,v)}=\big([A,B],\half(Av-Bu)\big).
\end{equation}
The quadruple $(\gl(V)\oplus
V,\pair{\cdot,\cdot},\Courant{\cdot,\cdot})$ is called the omni-Lie
algebra associated to the vector space $V$. A {\em Dirac structure}
of the omni-Lie algebra $\gl(V)\oplus V$ is a maximal isotropic
subspace on which the bracket $\Courant{\cdot, \cdot}$ becomes a Lie
bracket upon restriction. For any skew-symmetric bilinear operation
$[\cdot,\cdot]:V\wedge V\longrightarrow V$, the induced linear map
$\ad:V\longrightarrow\gl(V)$ is defined by
$$
\ad_u(v)=[u,v],\quad\forall~u,v\in V.
$$
The graph of the  map $\ad$, which we denote by $\frkG_\ad\subset
\gl(V)\oplus V$, is given by
\begin{equation}\label{eq:dirac}
 \frkG_\ad=\{(\ad_u,u)\in\gl(V)\oplus V|\forall~u\in V\}.
\end{equation}
Denote by $\id:\frkG_\ad\longrightarrow\gl(V)\oplus V$ the natural
embedding map. Obviously, $\frkG_\ad$ is a maximal isotropic
subspace of $\gl(V)\oplus V$ since the bilinear operation
$[\cdot,\cdot]$ is skew-symmetric. It is shown in
\cite{weinstein:omni} that $[\cdot,\cdot]$ is a Lie algebra
structure on $V$ if and only if $\frkG_\ad$ is a Dirac structure. In
this case, the map $\ad:V\longrightarrow \gl(V)$ is a Lie algebra
morphism. Consequently, the map
\begin{equation}\label{map:VadV}
V\to \frkG_\ad: v\mapsto (\ad_v,v),
\end{equation}
is a Lie algebra isomorphism.

 The factor of $\half$ in
(\ref{def:bracketomni}) spoils the Jacobi identity. More precisely,
we have
\begin{eqnarray*}
\Courant{\Courant{(A,u),(B,v)},(C,w)}+c.p.&=&\frac{1}{4}\big([A,B]w+[B,C]u+[C,A]v\big)\\
&\triangleq& T\big((A,u),(B,v),(C,w)\big).
\end{eqnarray*}
Thus $\Courant{\cdot,\cdot}$ is not a Lie bracket. However, the
Jacobiator is an exact term and we can extend the omni-Lie algebra
$\gl(V)\oplus V$ to the Lie 2-algebra whose degree-0 part is
$\gl(V)\oplus V$,
\begin{equation} \label{eq:omni-2}
\left\{\begin{array}{rlll}V& \stackrel{0+\Id}{\longrightarrow}&\gl(V)\oplus V,& \\
l_2(e_1,e_2)&=&\Courant{e_1,e_2},&\quad\mbox{for $e_1, e_2 \in \gl(V)\oplus V$},\\
l_2(e,f)&=&\Courant{e,\dM f},&\quad\mbox{for $e \in \gl(V)\oplus V, f \in V$},\\
l_3(e_1,e_2,e_3)&=&-T(e_1,e_2,e_3),&\quad\mbox{for $e_1, e_2, e_3
\in \gl(V) \oplus V$}.\end{array}\right.
\end{equation} such that the Jacobiator is measured by a ternary bracket taking value in the degree-1 part $V$.
This Lie 2-algebra is the semidirect product of the Lie algebra
$\gl(V)$ with the representation up to homotopy $(\mu, \nu)$
 on the 2-term
complex $V\stackrel{\Id}{\longrightarrow}V$,
$$
\mu(A)(u)=\half Au,\quad\nu(A,B)=\frac{1}{4}[A,B].
$$
Thus, we denote the Lie 2-algebra \eqref{eq:omni-2} by
$\gl(V)\ltimes(V\stackrel{\Id}{\longrightarrow}V)$.

$(\mu,\nu)$ is the morphism  from $\gl(V)$ to the crossed module
$(\gl(V), \gl(V), \Id, \ad)$.  The corresponding butterfly  is as
follows:
\begin{equation}
 \xymatrix{
0\ar[dd]\ar[dr]^{\kappa}&&\gl(V)\ar[dl]_{i_2}\ar[dd]\\
&\gl(V)\oplus\gl(V)\ar[dl]^{\sigma}\ar[dr]_{\rho}&\\
 \gl(V)&& \gl(V). }
\end{equation}
The Lie algebra structure on $\gl(V)\oplus\gl(V)$ is given by
\begin{equation}\label{bracket of omni}
[(A_1,A_2),(B_1,B_2)]=\big([A_1,B_1],\half([A_1,B_2]+[A_2,B_1])+[A_2,B_2]-\frac{1}{4}[A_1,B_1]\big),
\end{equation}
and the maps $\kappa,i_2,\sigma,\rho$ are given by
\begin{eqnarray*}
\kappa(0)&=&(0,0),\\
i_2(A)&=&(0,A),\\
\sigma(A_1,A_2)&=&A_1,\\
\rho(A_1,A_2)&=&\half A_1+A_2.
\end{eqnarray*}
By Theorem \ref{thm:equivalent}, we obtain a strict morphism
$(\psi_1,\psi_0)=(\Id,\rho)$ from the crossed module of Lie algebras
$(\gl(V),\gl(V)\oplus\gl(V),i_2,\ad)$ to $(\gl(V),\gl(V),\Id,\ad)$.
Furthermore, by Theorem \ref{thm:equivalent}, the semidirect product
\begin{equation} \label{eq:glv-glglv}
 (\gl(V)\ltimes V,(\gl(V)\oplus\gl(V))\ltimes
V,i_2\times\Id,\widetilde{\phi})
\end{equation}
 is equivalent to the Lie 2-algebra $\gl(V)\ltimes(V\stackrel{\Id}{\longrightarrow}V)$ via the morphism
$(f_0,f_1,f_2)$ given by
\begin{eqnarray}
\label{eq:omnif0}f_0(A,B,u)&=&(A,u),\\
\label{eq:omnif1}f_1(A,u)&=&u,\\
\nonumber f_2((A,B,u),(A^\prime,B^\prime,v))&=&Bv-B^\prime u.
\end{eqnarray}
Here the action $\widetilde{\phi}$ is given by
$$
\widetilde{\phi}_{(A,B,u)}(C,v)=\Big([\half A+B,C],(\half
A+B)(v)-Cu\Big),
$$
and the Lie algebra $\gl(V)\ltimes V$ is the semidirect product via
the natural action of $\gl(V)$ on $V$. The Lie algebra structure
on $(\gl(V)\oplus\gl(V))\ltimes V$ is given by
\begin{equation}\label{eq:bracket-glglv}
[(A_1,A_2,u),(B_1,B_2,v)]=\big([(A_1,A_2),(B_1,B_2)],(\half
A_1+A_2)(v)-(\half B_1+B_2)(u)\big).
\end{equation}

 By Corollary \ref{thm:integration}, we obtain the following
crossed module of Lie groups as the integration of the Lie 2-algebra
\eqref{eq:omni-2}
\begin{equation}\label{eqn:int omni}
(\mathcal{G}\ltimes V,\mathfrak{S}\ltimes V,(\int i_2)\times
\Id,\Phi),
\end{equation}
where $\mathcal{G}$ is the connected and simply connected Lie group
of the Lie algebra $\gl(V)$, $\mathfrak{S}$ is the connected and
simply connected Lie group of the Lie algebra $\gl(V)\oplus\gl(V)$
with the Lie bracket (\ref{bracket of omni}), $\Phi$ is the
integration of $\widetilde{\phi}$, and $\int i_2: \mathcal{G} \to
\mathfrak{S}$ is the unique map integrating $i_2: \gl(V) \to \gl(V)
\oplus \gl(V)$. \vspace{1mm}

We conclude by the following proposition.
\begin{pro}
The strict Lie 2-group corresponding to \eqref{eqn:int omni} is an
integration of the Lie 2-algebra
$\gl(V)\ltimes(V\stackrel{\Id}{\longrightarrow}V)$ associated to an
omni-Lie algebra $\gl(V)\oplus V$.
\end{pro}

\begin{rmk}In general, $GL(V)$ is neither connected nor simply connected.
It has two connected components $GL(V)_-$ and $GL(V)_+$ determined
by the sign of the determinant. When $\dim(V)=1$,
$\mathcal{G}=GL(V)_+=\R^{>0}$. When $\dim(V)=2$,
$\pi_1(GL(V)_+)=\mathbb Z$ and $\mathcal{G}$ is a $\mathbb Z$-cover
of $GL(V)_+$. When $\dim(V)\ge 3$, $\pi_1(GL(V)_+)=\mathbb Z_2$, and
$\mathcal{G}$ is a double cover of $GL(V)_+$.
\end{rmk}

\begin{rmk}
The Lie 2-algebra $\gl(V)\ltimes(V\stackrel{\Id}{\longrightarrow}V)$
is simply equivalent to the Lie algebra $\gl(V)$. Thus by our
definition, $GL(V)$ is also its integration. However, the omni-Lie
bracket $\Courant{\cdot, \cdot}$  does not appear in $\gl(V)$
anymore, while $\Courant{\cdot, \cdot}$ is a part of the structure
in the more complicated equivalent object \eqref{eq:glv-glglv} as
shown in \eqref{eq:bracket-glglv}. Moreover by Proposition
\ref{pro:sub-lie-2}, the Lie 2-algebra
$\gl(V)\ltimes(V\stackrel{\Id}{\longrightarrow}V)$ is a
sub-Lie-2-algebra of \eqref{eq:glv-glglv}. Thus we believe that
\eqref{eqn:int omni} is a more meaningful integration of
$\gl(V)\ltimes(V\stackrel{\Id}{\longrightarrow}V)$ than  simply
$GL(V)$.
\end{rmk}

We
must mention that much earlier than us, long before the program of
integration of $L_\infty$-algebras took shape in mathematics, Kinyon
and Weinstein set off to solve this problem in their pioneering work
\cite{kinyon-weinstein}. Their approach is to associate with the
omni-Lie algebra a Lie-Yamaguti algebra and use Kikkawa's
construction integrating Lie-Yamaguti algebras. They succeed in
finding a canonical left loop structure with non-associative
multiplication (which generalizes a group structure). However, as
the authors themselves point out, it is not satisfactory because it
does not obey the following testing requirement,
\begin{quote}
 \textit{Testing Requirement: the construction, when restricted to a
subspace where the Courant bracket becomes a Lie bracket, needs to
reproduce a group structure.}
               \end{quote}
They also hint at an alternative approach, in which one views the
Courant
bracket as a part of an $L_\infty$-structure (see Definition
\ref{def:Linfinity}),  as in \cite{rw}.  This is exactly the direction
we are taken now.

Now we study the testing requirement for our construction. By the
definition of Dirac structures,  we need to see that any Dirac
structure in the omni-Lie algebra $\gl(V)\oplus V$ integrates via
this procedure to a sub-Lie (1-)group of this Lie 2-group
\eqref{eqn:int omni}. As in \cite{kinyon-weinstein}, for simplicity,
we only consider the Dirac structure  $\frkG_\ad$ which is  a
graph as described in \eqref{eq:dirac}.

First of all, we view $\frkG_\ad$ as a trivial Lie 2-algebra
$0\stackrel{0}{\longrightarrow}{\frkG_\ad}$, which is a
sub-Lie-2-algebra of
$\gl(V)\ltimes(V\stackrel{\Id}{\longrightarrow}V)$. Our integration
of $\gl(V)\ltimes(V\stackrel{\Id}{\longrightarrow}V)$ consists first
of all in passing to an equivalent Lie 2-algebra which is strict,
and then in integrating the latter Lie 2-algebra. Since the
integration of equivalent Lie 2-algebras gives equivalent Lie
2-groups (see Section \ref{sec:butterfly}), the integration of our
Lie 2-algebra $0\stackrel{0}{\longrightarrow}{\frkG_\ad}$ is simply
$G(V)$, where $G(V)$ is the simply connected Lie group of the Lie
algebra $V$.

Now we make this more explicit by tracing the equivalence of Lie
2-algebras and give the precise equivalent sub-Lie-2-algebra. We
view the integration of the latter as a sub-Lie 2-group of the
integration \eqref{eqn:int omni}. Finally we verify this sub-Lie
2-group is equivalent to a Lie group. This part of calculation is
logically redundant, but we see more clearly some non-strict
phenomenon happening via pulling back a sub-Lie-2-algebra by an
equivalence.

The preimage under the map $(f_0,f_1)$ given by (\ref{eq:omnif0})
and (\ref{eq:omnif1}) of this sub-Lie-2-algebra
$0\stackrel{0}{\longrightarrow}{\frkG_\ad}$ is the sub-complex
$\gl(V)\stackrel{i_2}{\longrightarrow}(\ad_V\oplus\gl(V))\times V$,
where the vector space $(\ad_V\oplus\gl(V))\times V$ is given by
$$
(\ad_V\oplus\gl(V))\times V=\{(\ad_u,A,u)\in
(\gl(V)\oplus\gl(V))\ltimes V|~\forall~u\in V, A\in\gl(V)\},
$$
which is isomorphic to $\gl(V)\times V$.

It is not hard to see that the complex
$\gl(V)\stackrel{i_2}{\longrightarrow}(\ad_V\oplus\gl(V))\times V$
is not closed under the Lie brackets on $\gl(V)\ltimes V$ and
$(\gl(V)\oplus\gl(V))\ltimes V$ because we have
\begin{eqnarray*}
\nonumber&& [(\ad_u, A, u), (\ad_v, B, v)]\\&=&\Big(
\ad_{[u,v]},\half([\ad_u,B]+[A,\ad_v])+[A,B]+\frac{1}{4}\ad_{[u,v]},[u,v]+Av-Bu
\Big).
\end{eqnarray*}
However, this complex is the image of a strict Lie 2-algebra which
we define as follows. Consider the complex $\gl(V)\stackrel{
i_1}{\longrightarrow}\gl(V)\times V$, where  $  i_1(A)=(A,0)$, with
the following Lie bracket operation
\begin{eqnarray}
\label{bracket:glvw}~[A,A^\prime]&=&[A,A^\prime],\qquad\qquad\qquad\qquad\qquad\qquad\qquad\qquad\qquad\qquad \mbox{on $\gl(V)$},\\
\label{bracket:glvv}~[(A,u),(B,v)]&=&\Big(\half([\ad_u,B]+[A,\ad_v])+[A,B]+\frac{1}{4}\ad_{[u,v]},[u,v]\Big),
\quad\mbox{on $\gl(V)\times V$}.
\end{eqnarray}
Define the action $\phi$ of the Lie algebra $\gl(V)\times V$ on the
Lie algebra $\gl(V)$ by
$$
\phi_{(A,u)}(B)=([A,B]+\half[\ad_u,B], 0).
$$
Then we obtain a crossed module of Lie algebras:
\begin{equation}\label{crossed module wv}(\gl(V),\gl(V)\times V,i_1,\phi).\end{equation}
 Define
$\psi_0:\gl(V)\times V\longrightarrow(\gl(V)\oplus\gl(V))\times V$
by
$$
\psi_0(A,u)=(\ad_u,A,u)
$$
and  let $\psi_1:\gl(V)\longrightarrow\gl(V)\times V$ be  the
natural inclusion map. Furthermore, define
$\psi_2:\wedge^2(\gl(V)\times V)\longrightarrow\gl(V)\times V$ by
$$
\psi_2((A,u),(B,v))=(0, Av-Bu).
$$
Then it is not hard to see that $(\psi_0,\psi_1,\psi_2)$ is a Lie
2-algebra morphism.

To summarize, we have the following commutative diagram of Lie 2-algebras:
\begin{equation}\label{eq:pull-back-sub}
\xymatrix{ \gl(V) \to \gl(V)\times V \ar[r]^{(\psi_i)\qquad}
\ar[d]_{0, \ad_{\pr_V}\times \pr_V } & \gl(V) \ltimes V \to (\gl(V)
\oplus \gl(V) )
 \ltimes V \ar[d]_{\pr_V, \sigma \times \Id} \\
0\to \frkG_{\ad} \ar[r]_{(0,\id)}  & V\to\gl(V) \oplus V. }
\end{equation}
It is not hard to see that the vertical arrows are equivalences of
Lie 2-algebras.  Thus the crossed module of Lie algebras
\eqref{crossed module wv} is the pull-back of the sub-Lie-2-algebra
$0\to \frkG_{\ad}$ that we are interested in. The two Lie algebras
$\gl(V) $ and $\gl(V) \times V$ with the Lie brackets
(\ref{bracket:glvw}) and (\ref{bracket:glvv}) are both extensions
(the first a trivial one) of Lie algebras fitting in the following
diagram of Lie algebras
\[
 \xymatrix{0\ar[r]& \gl(V) \ar[d]_{\Id} \ar[r]^{\Id} & \gl(V) \ar[d]_{ i_1} \ar[r]^{\qquad0} &0 \ar[d]_{0} \ar[r] & 0\\
0\ar[r] &\gl(V) \ar[r]^{i_1} & \gl(V) \times V \ar[r]_{\qquad\pr_V}
& V \ar[r] & 0,}
\]
This diagram integrates to a commutative diagram of simply connected
Lie groups
\[
 \xymatrix{1\ar[r]& \huaG \ar[d]_{\Id} \ar[r]^{\Id} & \huaG   \ar[d]_{\int( i_1)} \ar[r]^{\quad1} & 1 \ar[d]_{1} \ar[r] & 1\\
1\ar[r] & \huaG \ar[r]^{} & G( \gl(V) \times V) \ar[r] & G(V) \ar[r]
& 1,}
\]
where $G( \gl(V) \times V)$ is the simply connected Lie group
integrating  the Lie algebra $\gl(V) \times V$ with  the Lie bracket
(\ref{bracket:glvv}). The Lie group $G( \gl(V) \times V)$ has  no
explicit form, however we  know that $\int (i_1)$ is injective since
$ i_1$ is injective. Thus by Corollary \ref{cor:equivalence-gp} the
integrated
 simply connected crossed module given by $\huaG  $ and $ G(\gl(V) \times V)$ is equivalent as a Lie 2-group to
the  Lie group $ G(\gl(V) \times V)/\huaG  \cong G(V)$.

In summary, rephrasing Kinyon-Weinstein's result in our higher language,  we can
associate a group-like object, which is a strict Lie 2-group, to an
omni-Lie algebra, meeting the testing requirement that Kinyon and
Weinstein ask. Moreover, given any such strict Lie 2-group coming
from integration, we can differentiate it to obtain a Lie 2-algebra
which contains the complete information of the omni-Lie algebra that
we start with. Thus it is appropriate to say that this strict Lie 2-group is the
integration of an omni-Lie algebra.

\bibliographystyle{habbrv}
\bibliography{../../bib/bibz}
\end{document}